\documentclass[preprint,12pt]{elsarticle}

\usepackage{amssymb}

\newtheorem{thm}{Theorem}

\newdefinition{dfn}{Definition}
\newdefinition{exm}{Example}
\newproof{pf}{Proof}

\newtheorem{lm}{Lemma}

\journal{Elsevier}

\begin{document}

\begin{frontmatter}



\title{Some New Paranormed  Difference Sequence Spaces Derived by Fibonacci Numbers}


\author[label1]{Emrah Evren Kara\corref{cor1}}
\ead{eevrenkara@hotmail.com}\cortext[cor1]{Corresponding Author
(Tel: +90 356 252 16 16, Fax: +90 356 252 15 85)}
\author[label2]{Serkan Demiriz}
\ead{serkandemiriz@gmail.com}

\address[label1]{Department of Mathematics, Duzce University,\\
 81620, Duzce, Turkey }
\address[label2]{Department of Mathematics,  Gaziosmanpa\c{s}a University,\\
 60250 Tokat, Turkey}

\begin{abstract}
In this study, we define new paranormed sequence spaces by the
sequences of Fibonacci numbers. Furthermore, we compute the
$\alpha-,\beta-$ and $\gamma-$ duals and obtain bases for these
sequence spaces. Besides this, we characterize the matrix
transformations from the new paranormed sequence spaces to the
Maddox's spaces $c_{0}(q),c(q),\ell(q)$ and $\ell_{\infty}(q)$.
\end{abstract}


\end{frontmatter}


\noindent
\section{Preliminaries,background and notation}
By $\omega$, we shall denote the space of all real valued sequences.
Any vector subspace of $\omega$ is called as a \textit{sequence
space}. We shall write $\ell_{\infty},c$ and $c_{0}$ for the spaces
of all bounded, convergent and null sequences, respectively. Also by
$bs,cs,\ell_{1}$ and $\ell_{p}$ ; we denote the spaces of all
bounded, convergent, absolutely and $p-$ absolutely convergent
series, respectively; $1< p< \infty$.

A linear topological space $X$ over the real field $\mathbb{R}$ is
said to be a paranormed space if there is a subadditive function
$g:X\rightarrow \mathbb{R}$ such that $g(\theta)=0, g(x)=g(-x)$ and
scalar multiplication is continuous,i.e.,
$|\alpha_{n}-\alpha|\rightarrow 0$ and $g(x_{n}-x)\rightarrow 0$
imply $g(\alpha_{n}x_{n}-\alpha x)\rightarrow 0$ for all $\alpha'$s
in $\mathbb{R}$ and all $x$'s in $X$, where $\theta$ is the zero
vector in the linear space $X$.


Assume here and after that  $(p_{k})$ be a bounded sequences of
strictly positive real numbers with $\sup p_{k}=H$ and $M=\max
\{1,H\}$. Then, the linear spaces $ c(p),c_{0}(p),\ell_{\infty}(p)$
and $\ell(p)$ were defined by Maddox \cite{m1,m2} (see also Simons
\cite{simons} and Nakano \cite{nakano}) as follows:

\begin{eqnarray*}
c(p)&=&\left\{x=(x_{k})\in\omega:\lim_{k\rightarrow \infty}
|x_{k}-l|^{p_{k}}=0 ~  \textrm{for some}~ l\in\mathbb{C}\right\},\\
c_{0}(p)&=&\left\{x=(x_{k})\in\omega:\lim_{k\rightarrow \infty}
|x_{k}|^{p_{k}}=0 \right\},\\
\ell_{\infty}(p)&=&\left\{x=(x_{k})\in\omega: \sup_{k\in
\mathbb{N}}|x_{k}|^{p_{k}}<\infty\right\}
\end{eqnarray*}
and
$$
\ell(p)=\bigg\{x=(x_{k})\in \omega: \sum_{k}
|x_{k}|^{p_{k}}<\infty\bigg\},
$$
which are the complete spaces paranormed by
$$
h_{1}(x)=\sup_{k\in\mathbb{N}} |x_{k}|^{p_{k}/M} \ {\rm iff} \
\inf_{p_{k}}>0 \qquad {\rm and} \qquad
h_{2}(x)=\bigg(\sum_{k}|x_{k}|^{p_{k}}\bigg)^{1/M},
$$
respectively. We shall assume throughout that
$p_{k}^{-1}+(p_{k}^{'})^{-1}=1$ provided $1<\inf p_{k}<H<\infty$.
For simplicity in notation, here and in what follows, the summation
without limits runs from $0$ to $\infty$. By $\mathcal{F}$ and
$\mathbb{N}_{k}$, we shall denote the collection of all finite
subsets of $\mathbb{N}$ and the set of all $n\in\mathbb{N}$ such
that $n\geq k$, respectively. We write by $\mathcal{U}$ for the set
of all sequences $u=(u_{n})$ such that $u_{n}\neq 0$ for all $n\in
\mathbb{N}$. For $u\in \mathcal{U}$, let $1/u=(1/u_{n})$.

For the sequence spaces $X$ and $Y$, define the set $S(X,Y)$ by
\begin{equation}\label{1.1}
S(X,Y)=\{z=(z_{k}): xz=(x_{k}z_{k})\in Y \   \ {\rm for\  \ all}\ \
x\in X\}.
\end{equation}
With the notation of (\ref{1.1}), the $\alpha-,\beta-$ and $\gamma-$
duals of a sequence space $X$, which are respectively denoted by
$X^{\alpha}, X^{\beta}$ and $X^{\gamma}$, are defined by
$$
X^{\alpha}=S(X,\ell_{1}) , \     \ X^{\beta}=S(X,cs) \ {\rm and} \ \
X^{\gamma}=S(X,bs).
$$

Let $(X,h)$ be a paranormed space. A sequence $(b_{k})$ of the
elements of $X$ is called a basis for $X$ if and only if, for each
$x\in X$, there exists a unique sequence $(\alpha_{k})$ of scalars
such that
$$
h\left (x-\sum_{k=0}^{n} \alpha_{k}b_{k}\right) \rightarrow 0 \ \ as
\     \   n\rightarrow\infty.
$$
The series $\sum \alpha_{k} b_{k}$ which has the sum $x$ is then
called the expansion of $x$ with respect to $(b_{n})$ and written as
$x=\sum \alpha_{k} b_{k}$.
Let $X,Y$ be any two sequence spaces and $A=(a_{nk})$ be an infinite
matrix of real numbers $a_{nk}$,where $n,k\in \mathbb{N}$. Then, we
say that $A$ defines a matrix mapping from $X$ into $Y$, and we
denote it by writing $A:X\rightarrow Y$, if for every sequence
$x=(x_{k})\in X$ the sequence $Ax=((Ax)_{n})$, the $A$-transform of
$x$, is in $Y$, where
\begin{equation}\label{1.2}
(Ax)_{n}=\sum_{k} a_{nk}x_{k},\     \ (n\in\mathbb{N}).
\end{equation}
By $(X:Y)$, we denote the class of all matrices $A$ such that
$A:X\rightarrow Y$. Thus, $A\in(X:Y)$ if and only if the series on
the right-hand side of (\ref{1.2}) converges for each $n\in
\mathbb{N}$ and every $x\in X$, and we have $Ax=\{(Ax)_{n}\}_{n\in
\mathbb{N}}\in Y$ for all $x\in X$. A sequence $x$ is said to be
$A$- summable to $\alpha$ if $Ax$ converges to $\alpha$ which is
called as the $A$- limit of $x$.



For a sequence space $X$, the matrix domain $X_{A}$ of an infinite
matrix $A$ is defined by
$$
X_{A}=\{x=(x_{k})\in\omega:\  \ Ax\in X\}.
$$
The approach constructing a new paranormed sequence space by means
of the matrix domain of a particular limitation method has recently
been employed by Malkowsky \cite{em}, Altay and Ba\c{s}ar
\cite{bafb3,bafb4}, F. Ba\c{s}ar et al., \cite{fbbamm}, Ayd{\i}n and
Ba\c{s}ar \cite{cafb1,cafb2}.

Define the sequence $\{f_n\}_{n=0}^{\infty}$ of Fibonacci numbers
given by the linear recurrence relations
$$
f_{0}=f_{1}=1 \ \textrm{and} \ f_{n}=f_{n-1}+f_{n-2}, \quad n\geq 2.
$$
In modern science and particularly physics, there is quite an
interest in the theory and applications of Fibonacci numbers . The
ratio of the successive Fibonacci numbers is as known golden ratio.
There are many applications of golden ratio in many places of
mathematics and physics, in general theory of high energy particle
theory \cite{fibo2}. Also, some basic properties of Fibonacci
numbers \cite{fibo1} are given as follows:
$$
\lim_{n\rightarrow \infty}
\frac{f_{n+1}}{f_{n}}=\frac{1+\sqrt{5}}{2}=\alpha \qquad
(\textrm{golden ratio})
$$
$$
\sum_{k=0}^{n} f_{k}=f_{n+2}-1 \qquad (n\in \mathbb{N}) \quad
\textrm{and} \quad \sum_{k} \frac{1}{f_{k}} \quad \textrm{converges}
$$
$$
f_{n-1}f_{n+1}-f_{n}^{2}=(-1)^{n+1} \quad (n\geq 1) \quad
(\textrm{Cassini formula}).
$$
Substituting for $f_{n+1}$ in Cassini's formula yields
$f_{n-1}^{2}+f_{n}f_{n-1}-f_{n}^{2}=(-1)^{n+1}$.


Let $f_n$ be the $n$th Fibonacci number for every $n\in \mathbb{N}$.
Then, the infinite Fibonacci matrix $\widehat{F}=(\widehat{f}_{nk})$
is defined by
$$
\widehat{f}_{nk}=\left\{\begin{array}{ll}
  \displaystyle -\frac{f_{n+1}}{f_{n}} & (k=n-1),\\
  \displaystyle \frac{f_{n}}{f_{n+1}} & (k=n),\\
  \displaystyle 0 & (0\leq k<n-1 \ \textrm{or} \ k>n)
\end{array}\right.
$$
where $n,k\in \mathbb{N}$ \cite{eek}.


The main purpose of this study is to introduce the  sequence spaces
$c_{0}(\widehat{F},p),c(\widehat{F},p)$
,$\ell_{\infty}(\widehat{F},p)$ and $\ell(\widehat{F},p)$ which are
the set of all sequences whose $\widehat{F}-$transforms are in the
spaces $c_{0}(p),c(p),\ell_{\infty}(p)$ and $\ell(p)$, respectively.
Also, we have investigated some topological structures, which have
completeness, the $\alpha-,\beta-$ and $\gamma-$ duals, and the
bases of these sequence spaces. Besides this, we characterize some
matrix mappings on these spaces.
\section{The Paranormed Fibonacci Difference Sequence Spaces }

In this section, we define the new  sequence spaces
$c_{0}(\widehat{F},p),c(\widehat{F},p)$
,$\ell_{\infty}(\widehat{F},p)$ and $\ell(\widehat{F},p)$  by using
the sequences of Fibonacci numbers, and prove that these sequence
spaces are the complete paranormed linear metric spaces and compute
their $\alpha-,\beta-$ and $\gamma-$ duals. Moreover, we give the
basis for the spaces $c_{0}(\widehat{F},p),c(\widehat{F},p)$ and
$\ell(\widehat{F},p)$ .

For a sequence space $X$, the matrix domain $X_{A}$ of an infinite
matrix $A$ is defined by
\begin{equation}\label{2.1}
X_{A}=\{x=(x_{k})\in\omega:\  \ Ax\in X\}.
\end{equation}

In \cite{cm}, Choudhary and Mishra have defined the sequence space
$\overline{\ell(p)}$ which consists of all sequences such that
$S$-transforms are in $\ell(p)$, where $S=(s_{nk})$ is defined by
$$
s_{nk}=\left\{\begin{array}{cc}
  \displaystyle 1, & (0 \leq k \leq n), \\
  0, & (k>n).
\end{array}\right.
$$
Ba\c{s}ar and Altay \cite{fbba1} have recently examined the space
$bs(p)$ which is formerly defined by Ba\c{s}ar in \cite{fb} as the
set of all series whose sequences of partial sums are in
$\ell_{\infty}(p)$. More recently, Altay and Ba\c{s}ar have studied
the sequence spaces $r^{t}(p),r_{\infty}^{t}(p)$ in \cite{bafb1} and
$r_{c}^{t}(p),r_{0}^{t}(p)$ in \cite{bafb2} which are derived by the
Riesz means from the sequence spaces $\ell(p),\ell_{\infty}(p),c(p)$
and $c_{0}(p)$ of Maddox, respectively. With the notation of
(\ref{2.1}), the spaces
$\overline{\ell(p)},bs(p),r^{t}(p),r_{\infty}^{t}(p),r_{c}^{t}(p)$
and $r_{0}^{t}(p)$ may be redefined by

$$
\overline{\ell(p)}=[\ell(p)]_{S},\   \
bs(p)=[\ell_{\infty}(p)]_{S},\   \ r^{t}(p)=[\ell(p)]_{R^{t}},
$$

$$
r_{\infty}^{t}(p)=[\ell_{\infty}(p)]_{R^{t}},\   \
r_{c}^{t}(p)=[c(p)]_{R^{t}},\   \ r_{0}^{t}(p)=[c_{0}(p)]_{R^{t}}.\\
$$


Following Choudhary and Mishra \cite{cm}, Ba\c{s}ar and Altay
\cite{fbba1}, Altay and Ba\c{s}ar \cite{bafb1,bafb2}, we define the
sequence spaces $c_{0}(\widehat{F},p),c(\widehat{F},p)$
,$\ell_{\infty}(\widehat{F},p)$ and $\ell(\widehat{F},p)$  by
\begin{eqnarray*}
c_{0}(\widehat{F},p)&=&\bigg\{x=(x_{k})\in \omega:
\lim_{n\rightarrow \infty}
\bigg|\frac{f_{n}}{f_{n+1}}x_{n}-\frac{f_{n+1}}{f_{n}}x_{n-1}\bigg|^{p_n}=0\bigg\}\\
c(\widehat{F},p)&=&\bigg\{x=(x_{k})\in \omega: \exists l\in
\mathbb{C} \ni \lim_{n\rightarrow \infty}
\bigg|\frac{f_{n}}{f_{n+1}}x_{n}-\frac{f_{n+1}}{f_{n}}x_{n-1}-l\bigg|^{p_n}=0\bigg\}\\
\ell_{\infty}(\widehat{F},p)&=&\bigg\{x=(x_{k})\in \omega:
\sup_{n\in \mathbb{N}}
\bigg|\frac{f_{n}}{f_{n+1}}x_{n}-\frac{f_{n+1}}{f_{n}}x_{n-1}\bigg|^{p_n}<\infty
\bigg\}
\end{eqnarray*}
and
$$
\ell(\widehat{F},p)=\bigg\{x=(x_{k})\in \omega: \sum_{n}
\bigg|\frac{f_{n}}{f_{n+1}}x_{n}-\frac{f_{n+1}}{f_{n}}x_{n-1}\bigg|^{p_n}<\infty\bigg\}.
$$
In the case $(p_{k})=e=(1,1,1,...),$ the sequence spaces
$c_{0}(\widehat{F},p),c(\widehat{F},p)$
,$\ell_{\infty}(\widehat{F},p)$ and $\ell(\widehat{F},p)$ are ,
respectively, reduced to the sequence spaces
$c_{0}(\widehat{F}),c(\widehat{F}),\ell_{\infty}(\widehat{F})$ and
$\ell_{p}(\widehat{F})$ which are introduced by E.E.Kara \cite{eek}
and M. Ba\c{s}ar{\i}r et al. \cite{eek1}.

With the notation (\ref{2.1}), we may redefine the spaces
$c_{0}(\widehat{F},p),c(\widehat{F},p)$
,$\ell_{\infty}(\widehat{F},p)$ and $\ell(\widehat{F},p)$ as
follows:
$$
c_{0}(\widehat{F},p)=\{c_{0}(p)\}_{\widehat{F}}, \quad
c(\widehat{F},p)=\{c(p)\}_{\widehat{F}},
$$
$$
\ell_{\infty}(\widehat{F},p)=\{\ell_{\infty}(p)\}_{\widehat{F}},
\quad \ell(\widehat{F},p)=\{\ell(p)\}_{\widehat{F}}.
$$

Define the sequence $y=(y_{k})$, which will be frequently used as
the $\widehat{F}-$transform of a sequence $x=(x_{k})$, i.e.
\begin{equation}\label{2.2}
y_{k}=\widehat{F}_{k}(x)=\frac{f_{k}}{f_{k+1}}x_{k}-\frac{f_{k+1}}{f_{k}}x_{k-1};
\quad (k\in \mathbb{N}_{0}).
\end{equation}
Since the proof may also be obtained in the similar way as for the
other spaces, to avoid the repetition of the similar statements, we
give the proof only for one of those spaces. Now, we may begin with
the following theorem which is essential in the study.


\begin{thm}\label{t2.1}
(i) The sequence spaces $c_{0}(\widehat{F},p),c(\widehat{F},p)$ and
$\ell_{\infty}(\widehat{F},p)$ are the complete linear metric spaces
paranormed by $g$, defined by
$$
g(x)=\sup_{k\in \mathbb{N}}
\bigg|\frac{f_{k}}{f_{k+1}}x_{k}-\frac{f_{k+1}}{f_{k}}x_{k-1}
\bigg|^{p_{k}/M}.
$$
$g$ is a paranorm for the spaces $c(\widehat{F},p)$ and
$\ell_{\infty}(\widehat{F},p)$ only in the trivial case $\inf
p_{k}>0$ when $c(\widehat{F},p)=c(\widehat{F})$ and
$\ell_{\infty}(\widehat{F},p)=\ell_{\infty}(\widehat{F})$.\\
(ii) $\ell_{p}(\widehat{F})$ is a complete linear metric space
paranormed by
$$
g^{*}(x)=\bigg(\sum_{k}
\bigg|\frac{f_{k}}{f_{k+1}}x_{k}-\frac{f_{k+1}}{f_{k}}x_{k-1}
\bigg|^{p_{k}}\bigg)^{1/M}.
$$
\end{thm}


\begin{pf}
We prove the theorem for the space $c_{0}(\widehat{F},p)$. The
linearity of $c_{0}(\widehat{F},p)$ with respect to the
coordinatewise addition and scalar multiplication follows from the
following inequalities which are satisfied for $x,z\in
c_{0}(\widehat{F},p)$ (see \cite [p.30]{m}):
\begin{eqnarray}\label{2.3}
\sup_{k\in
\mathbb{N}}\bigg|\frac{f_{k}}{f_{k+1}}(x_{k}+z_{k})-\frac{f_{k+1}}{f_{k}}(x_{k-1}+z_{k-1})
\bigg|^{p_{k}/M} \nonumber &\leq&\sup_{k\in
\mathbb{N}}\bigg|\frac{f_{k}}{f_{k+1}}x_{k}-\frac{f_{k+1}}{f_{k}}x_{k-1}
\bigg|^{p_{k}/M}\\
&+&\sup_{k\in
\mathbb{N}}\bigg|\frac{f_{k}}{f_{k+1}}z_{k}-\frac{f_{k+1}}{f_{k}}z_{k-1}
\bigg|^{p_{k}/M}\nonumber\\
\end{eqnarray}
and for any $\alpha \in \mathbb{R}$ (see \cite{m2}),
\begin{equation}\label{2.4}
|\alpha|^{p_{k}}\leq \max \{1,|\alpha|^{M}\}.
\end{equation}
It is clear that $g(\theta)=0$ and $g(x)=g(-x)$ for all $x\in
c_{0}(\widehat{F},p)$. Again the inequalities (\ref{2.3}) and
(\ref{2.4}) yield the subadditivity of $g$ and
$$
g(\alpha x)\leq \max \{1,|\alpha|\}g(x).
$$
Let $\{x^{n}\}$ be any sequence of the points $x^{n} \in
c_{0}(\widehat{F},p)$ such that $g(x^{n}-x)\rightarrow 0$ and
$(\alpha_{n})$ also be any sequence of scalars such that
$\alpha_{n}\rightarrow \alpha$. Then, since the inequality
$$
g(x^{n})\leq g(x)+g(x^{n}-x)
$$
holds by the subadditivity of $g$, $\{g(x^{n})\}$ is bounded and we
thus have
\begin{eqnarray*}
g(\alpha_{n} x^{n}-\alpha x)&=&
\sup_{k\in \mathbb{N}} \bigg|\frac{f_k}{f_{k+1}}(\alpha_{n}x_{k}^{n}-\alpha x_{k})-\frac{f_{k+1}}{f_{k}}(\alpha_{n}x_{k-1}^{n}-\alpha x_{k-1})\bigg|^{p_{k}/M}\\
&\leq& |\alpha_{n}-\alpha|\ \ g(x^{n})+ |\alpha|\ \ g(x^{n}-x),
\end{eqnarray*}
which tends to zero as $n\rightarrow \infty$. That is to say that
the scalar multiplication is continuous. Hence, $g$ is a paranorm on
the space $c_{0}(\widehat{F},p)$.

It remains to prove the completeness of the space $
c_{0}(\widehat{F},p)$. Let $\{x^{i}\}$ be any Cauchy sequence in the
space $c_{0}(\widehat{F},p)$, where $x^{i}=\{x_{0}^{(i)},
x_{1}^{(i)}, ...\}$. Then, for a given $\varepsilon>0$ there exists
a positive integer $n_{0}(\varepsilon)$ such that
$$
g(x^{i}-x^{j})< \frac{\varepsilon}{2}
$$
for all $i,j\geq n_{0}(\varepsilon)$. We obtain by using definition
of $g$ for each fixed $k\in \mathbb{N}$ that
\begin{eqnarray}\label{2.5}
\big|\{\widehat{F}x^{i}\}_{k}-\{\widehat{F}x^{j}\}_{k}\big|^{p_{k}/M}
&\leq& \sup_{k\in \mathbb{N}}
\big|\{\widehat{F}x^{i}\}_{k}-\{\widehat{F}x^{j}\}_{k}\big|^{p_{k}/M}\nonumber\\
&<&\frac{\varepsilon}{2}
\end{eqnarray}
for every $i,j\geq n_{0}(\varepsilon)$, which leads us to the fact
that $\{(\widehat{F}x^{0})_{k},$ $(\widehat{F}x^{1})_{k},...\}$ is a
Cauchy sequence of real numbers for every fixed $k\in \mathbb{N}$.
Since $\mathbb{R}$ is complete, it converges, say
$$
\{\widehat{F}x^{i}\}_{k}\rightarrow \{\widehat{F}x\}_{k}
$$
as $i\rightarrow \infty$. Using these infinitely many limits
$(\widehat{F}x)_{0},(\widehat{F}x)_{1},...$, we define the sequence
$\{(\widehat{F}x)_{0},(\widehat{F}x)_{1},...\}$. We have from
(\ref{2.5}) with $j\rightarrow \infty$ that
\begin{equation}\label{2.6}
\big|\{\widehat{F}x^{i}\}_{k}-\{\widehat{F}x\}_{k}\big|^{p_{k}/M}\leq
\frac{\varepsilon}{2} \quad (i\geq n_{0}(\varepsilon))
\end{equation}
for every fixed $k\in \mathbb{N}$. Since $x^{i}=\{x_{k}^{(i)}\}\in
c_{0}(\widehat{F},p)$,
$$
\big|\{\widehat{F}x^{i}\}_{k}\big|^{p_{k}/M}<\frac{\varepsilon}{2}
$$
for all $k\in \mathbb{N}$. Therefore, we obtain (\ref{2.6}) that
\begin{eqnarray*}
\big|\{\widehat{F}x\}_{k}\big|^{p_{k}/M}&\leq&
\big|\{\widehat{F}x\}_{k}-\{\widehat{F}x^{i}\}_{k}\big|^{p_{k}/M}
+\big|\{\widehat{F}x^{i}\}_{k}\big|^{p_{k}/M}\\
&<&\varepsilon  \quad (i\geq n_{0}(\varepsilon)).
\end{eqnarray*}
This shows that the sequence $\{\widehat{F}x\}$ belongs to  the
space $c_{0}(p)$. Since $\{x^{i}\}$ was an arbitrary Cauchy
sequence, the space $c_{0}(\widehat{F},p)$ is complete and this
concludes the proof.
\end{pf}

Therefore, one can easily check that the absolute property does not
hold on the spaces $c_{0}(\widehat{F},p),c(\widehat{F},p)$ ,
$\ell_{\infty}(\widehat{F},p)$ and $\ell(\widehat{F},p)$ that is
$h(x)\neq h(|x|)$ for at least one sequence in those spaces, and
this says that $c_{0}(\widehat{F},p),c(\widehat{F},p)$
,$\ell_{\infty}(\widehat{F},p)$ and $\ell(\widehat{F},p)$ are the
sequence spaces of non-absolute type; where $|x|=(|x_k|)$.

\begin{thm}\label{t2.2}
The sequence spaces $c_{0}(\widehat{F},p),c(\widehat{F},p)$
,$\ell_{\infty}(\widehat{F},p)$ and $\ell(\widehat{F},p)$ are
linearly isomorphic to the spaces $c_{0}(p),c(p),\ell_{\infty}(p)$
and $\ell(p)$, respectively, where $0<p_{k}\leq H<\infty$.
\end{thm}
\begin{pf}
We establish this for the space $\ell_{\infty}(\widehat{F},p)$. To
prove the theorem, we should show the existence of a linear
bijection between the spaces $\ell_{\infty}(\widehat{F},p)$ and
$\ell_{\infty}(p)$ for $0<p_{k}\leq H<\infty$. With the notation of
(\ref{2.2}), define the transformations $T$ from
$\ell_{\infty}(\widehat{F},p)$ to $\ell_{\infty}(p)$ by $x\mapsto
y=Tx$. The linearity of $T$ is trivial. Further, it is obvious that
$x=\theta$ whenever $Tx=\theta$ and hence $T$ is injective.

Let $y=(y_{k})\in \ell_{\infty}(p)$ and define the sequence
$x=(x_{k})$ by
$$
x_{k}=\sum_{j=0}^{k} \frac{f_{k+1}^{2}}{f_{j}f_{j+1}}y_{j}; \quad
(k\in \mathbb{N}).
$$
Then, we get that
\begin{eqnarray*}
g(x)&=&\sup_{k\in
\mathbb{N}}\bigg|\frac{f_{k}}{f_{k+1}}x_{k}-\frac{f_{k+1}}{f_{k}}x_{k-1}
\bigg|^{p_{k}/M}\\
&=&\sup_{k\in \mathbb{N}}\bigg|\frac{f_{k}}{f_{k+1}}\sum_{j=0}^{k}
\frac{f_{k+1}^{2}}{f_{j}f_{j+1}}y_{j}-\frac{f_{k+1}}{f_{k}}\sum_{j=0}^{k-1}
\frac{f_{k}^{2}}{f_{j}f_{j+1}}y_{j} \bigg|^{p_{k}/M}\\
&=& \sup_{k\in \mathbb{N}}|y_{k}|^{p_{k}/M}=h_{1}(y)<\infty.
\end{eqnarray*}
Thus, we deduce that $x\in \ell_{\infty}(\widehat{F},p)$ and
consequently $T$ is surjective and is paranorm preserving. Hence,
$T$ is a linear bijection and this says us that the spaces
$\ell_{\infty}(\widehat{F},p)$ and $\ell_{\infty}(p)$ are linearly
isomorphic, as desired.
\end{pf}

We shall quote some lemmas which are needed in proving related to
the duals our theorems.
\begin{lm}\label{l2.3}\cite [Theorem 5.1.1 with $q_n=1$]{kgge}
$A\in (c_{0}(p):\ell(q))$  if and only if
\begin{equation}
\sup_{K \in \mathcal{F}} \sum_{n}\left|\sum_{k\in K}
a_{nk}B^{-1/p_{k}}\right|<\infty, \    \ (\exists B\in
\mathbb{N}_{2}).
\end{equation}
\end{lm}
\begin{lm}\label{l2.4}\cite [Theorem 5.1.9 with $q_n=1$]{kgge}
$A\in (c_{0}(p):c(q))$  if and only if
\begin{equation}
\sup_{n \in \mathcal{\mathbb{N}}} \sum_{k}
|a_{nk}|B^{-1/p_{k}}<\infty \ \ (\exists B\in \mathbb{N}_{2}),
\end{equation}
\begin{equation}
\exists (\alpha_{k})\subset \mathbb{R}\ni \lim _{n\rightarrow
\infty} |a_{nk}-\alpha_{k}|=0\    \ for\  \ all\  k\in \mathbb{N},
\end{equation}
\begin{equation}
\exists (\alpha_{k})\subset \mathbb{R}\ni \sup _{n\in \mathbb{N}}
\sum_{k}|a_{nk}-\alpha_{k}|B^{-1/p_{k}}<\infty. \quad (\exists B\in
\mathbb{N}_{2})
\end{equation}
\end{lm}

\begin{lm}\label{l2.5}\cite [Theorem 5.1.13 with $q_n=1$]{kgge}
$A\in (c_{0}(p):\ell_{\infty}(q))$  if and only if
\begin{equation}
\sup_{n\in \mathbb{N}}\sum_{k} |a_{nk}|B^{-1/p_{k}}<\infty. \quad
(\exists B\in \mathbb{N}_{2})
\end{equation}
\end{lm}

\begin{lm} \label{l2.6}\cite [Theorem 5.1.0 with $q_{n}=1$]{kgge}(i) Let $1<p_{k}\leq
H<\infty$ for all $k\in \mathbb{N}$. Then, $A\in (\ell(p):\ell_{1})$
if and only if there exists an integer $B>1$ such that
\begin{equation}\label{2.11}
\sup_{K \in \mathcal{F}} \sum_{k}\left|\sum_{n\in K}
a_{nk}B^{-1}\right|^{p_{k}^{'}}<\infty .
\end{equation}
(ii) Let $0<p_{k}\leq 1$ for all $k\in \mathbb{N}$. Then, $A\in
(\ell(p):\ell_{1})$ if and only if
\begin{equation}\label{2.12}
\sup_{K \in \mathcal{F}} \sup_{k\in \mathbb{N}}\left|\sum_{n\in K}
a_{nk}\right|^{p_{k}}<\infty .
\end{equation}
\end{lm}
\begin{lm} \label{l2.7}\cite [ Theorem 1 (i)-(ii)]{kgge} (i) Let $1<p_{k}\leq
H<\infty$ for all $k\in \mathbb{N}$. Then, $A\in
(\ell(p):\ell_{\infty})$ if and only if there exists an integer
$B>1$ such that

\begin{equation}\label{2.13}
\sup_{n \in \mathbb{N} } \sum_{k} |a_{nk}B^{-1}|^{p_{k}^{'}}<\infty
.
\end{equation}
(ii) Let $0<p_{k}\leq 1$ for all $k\in \mathbb{N}$. Then, $A\in
(\ell(p):\ell_{\infty})$ if and only if
\begin{equation}\label{2.14}
\sup_{n,k \in \mathbb{N}}|a_{nk}|^{p_{k}}<\infty.
\end{equation}
\end{lm}
\begin{lm}\label{l2.8}\cite [Corollary for Theorem 1]{kgge} Let $0<p_{k}\leq H<\infty$ for all $k\in
\mathbb{N}$. Then,  $A\in (\ell(p):c)$  if and only if (\ref{2.13}),
(\ref{2.14}) hold, and
\begin{equation}\label{2.15}
\lim_{n\rightarrow \infty} a_{nk}=\beta_{k}, \quad (k\in \mathbb{N})
\end{equation}
also holds.
\end{lm}






\begin{thm}\label{t2.9}
Let $K^{*}=\{k\in \mathbb{N}: 0\leq k\leq n\}\cap K$ for $K\in
\mathcal{F}$ and $B\in \mathbb{N}_{2}$. Define the sets
$\widehat{F}_{1}(p),\widehat{F}_{2}(p),\widehat{F}_{3}(p),\widehat{F}_{4}(p),\widehat{F}_{5}(p),\widehat{F}_{6}(p),\widehat{F}_{7}(p)$
and $\widehat{F}_{8}(p)$ as follows:
\begin{eqnarray*}
\widehat{F}_{1}(p)&=&\bigcup_{B>1} \bigg\{a=(a_{k})\in \omega:
\sup_{K\in \mathcal{F}} \sum_{n} \bigg|\sum_{k\in K^{*}}
\frac{f_{n+1}^{2}}{f_{k}f_{k+1}}a_{n}B^{-1/p_{k}}\bigg|<\infty
\bigg\}\\
\widehat{F}_{2}(p)&=&\bigg\{a=(a_{k})\in \omega: \sum_{n}
\bigg|\sum_{k=0}^{n}
\frac{f_{n+1}^{2}}{f_{k}f_{k+1}}a_{n}\bigg|<\infty\bigg\}\\
\widehat{F}_{3}(p)&=&\bigcup_{B>1} \bigg\{a=(a_{k})\in \omega:
\sup_{n\in \mathbb{N}} \sum_{k=0}^{n}\bigg|\sum_{j=k}^{n}
\frac{f_{j+1}^{2}}{f_{k}f_{k+1}}a_{j}\bigg|B^{-1/p_{k}}<\infty
\bigg\}\\
\widehat{F}_{4}(p)&=&\bigg\{a=(a_{k})\in \omega:
\bigg|\sum_{j=k}^{\infty}
\frac{f_{j+1}^{2}}{f_{k}f_{k+1}}a_{j}\bigg|<\infty \quad
\textrm{for all} \quad k\in \mathbb{N}\bigg\}\\
\widehat{F}_{5}(p)&=& \bigcup_{B>1} \bigg\{a=(a_{k})\in \omega:
\exists (\alpha_{k})\subset \mathbb{R}\ni \sup_{n\in \mathbb{N}}
\sum_{k=0}^{n} \bigg|\sum_{j=k}^{n}
\frac{f_{j+1}^{2}}{f_{k}f_{k+1}}a_{j}-\alpha_{k}\bigg|B^{-1/p_{k}}<\infty
\bigg\}\\
\widehat{F}_{6}(p)&=&\bigg\{a=(a_{k})\in \omega: \exists \alpha \in
\mathbb{R}\ni \lim_{n\rightarrow \infty} \bigg|\sum_{k=0}^{n}
\sum_{j=k}^{n} \frac{f_{j+1}^{2}}{f_{k}f_{k+1}}a_{j}-\alpha
\bigg|=0\bigg\}\\
\widehat{F}_{7}(p)&=&\bigg\{a=(a_{k})\in \omega: \sup_{n\in
\mathbb{N}} \bigg|\sum_{k=0}^{n} \sum_{j=k}^{n}
\frac{f_{j+1}^{2}}{f_{k}f_{k+1}}a_{j}\bigg|<\infty\bigg\}
\end{eqnarray*}
Then,

(i) $\{c_{0}(\widehat{F},p)\}^{\alpha}=\widehat{F}_{1}(p)$ \qquad
(ii) $\{c(\widehat{F},p)\}^{\alpha}=\widehat{F}_{1}(p)\cap
\widehat{F}_{2}(p)$ \\

(iii) $\{c_{0}(\widehat{F},p)\}^{\beta}=\widehat{F}_{3}(p)\cap
\widehat{F}_{4}(p)\cap \widehat{F}_{5}(p)$ \\

(iv) $\{c(\widehat{F},p)\}^{\beta}=
\{c_{0}(\widehat{F},p)\}^{\beta}\cap \widehat{F}_{6}(p)$\\

(v) $\{c_{0}(\widehat{F},p)\}^{\gamma}=\widehat{F}_{3}(p)$ \quad
(vi) $\{c(\widehat{F},p)\}^{\gamma}=\widehat{F}_{3}(p)\cap
\widehat{F}_{7}(p)$
\end{thm}
\begin{pf}
We give the proof for the space $c_{0}(\widehat{F},p)$. Let us take
any $a=(a_{n})\in \omega$ and define the matrix $C=(c_{nk})$ via the
sequence $a=(a_{n})$ by
$$
c_{nk}=\left\{\begin{array}{ll}
  \displaystyle \frac{f_{n+1}^{2}}{f_{k}f_{k+1}}a_{n}, & 0 \leq k \leq
  n,\\
   0, & k>n,
\end{array}\right.
$$
where $n,k\in \mathbb{N}$. Bearing in mind (\ref{2.2}) we
immediately derive that
\begin{eqnarray}\label{2.16}
a_{n}x_{n}=\sum_{k=0}^{n}
\frac{f_{n+1}^{2}}{f_{k}f_{k+1}}a_{n}y_{k}=(Cy)_{n}; \quad (n\in
\mathbb{N}).
\end{eqnarray}
We therefore observe by (\ref{2.16}) that $ax=(a_{n}x_{n})\in
\ell_{1}$ whenever $x\in c_{0}(\widehat{F},p)$ if and only if $Cy\in
\ell_{1}$ whenever $y\in c_{0}(p)$. This means that $a=(a_{n})\in
\{c_{0}(\widehat{F},p)\}^{\alpha}$ whenever $x=(x_{n})\in
c_{0}(\widehat{F},p)$ if and only if $C\in (c_{0}(p):\ell_{1})$.
Then, we derive by Lemma \ref{l2.3} that
$$
\{c_{0}(\widehat{F},p)\}^{\alpha}=\widehat{F}_{1}(p).
$$

Consider the equation for $n\in \mathbb{N}$,
\begin{eqnarray}\label{2.17}
\sum_{k=0}^{n}a_{k}x_{k}&=&\sum_{k=0}^{n}a_{k}
\bigg(\sum_{j=0}^{n}\frac{f_{k+1}^{2}}{f_{j}f_{j+1}}y_{j}\bigg)\nonumber\\
&=&\sum_{k=0}^{n} \bigg(\sum_{j=k}^{n}
\frac{f_{j+1}^{2}}{f_{k}f_{k+1}}a_{j}\bigg)y_{k} \nonumber\\
&=&(Dy)_{n}
\end{eqnarray}
where $D=(d_{nk})$ is defined by
$$
d_{nk}=\left\{\begin{array}{ll}
  \displaystyle \sum_{j=k}^{n}
\frac{f_{j+1}^{2}}{f_{k}f_{k+1}}a_{j}, & 0 \leq k \leq n,
  \\
    0, & k>n,
\end{array}\right.
$$
where $n,k\in \mathbb{N}$. Thus, we deduce from Lemma \ref{l2.4}
with (\ref{2.17}) that $ax=(a_{k}x_{k})\in cs$ whenever
$x=(x_{k})\in c_{0}(\widehat{F},p)$ if and only if $D y\in c$
whenever $y=(y_{k})\in c_{0}(p)$. This means that $a=(a_{n})\in
\{c_{0}(\widehat{F},p)\}^{\beta}$ whenever $x=(x_{n})\in
c_{0}(\widehat{F},p)$ if and only if $D\in (c_{0}(p):c)$. Therefore
we derive from Lemma \ref{l2.4} that
$$\{c_{0}(\widehat{F},p)\}^{\beta}=\widehat{F}_{3}(p)\cap
\widehat{F}_{4}(p)\cap \widehat{F}_{5}(p).$$

As this, we deduce from Lemma \ref{l2.5} with (\ref{2.17}) that
$ax=(a_{k}x_{k})\in bs$ whenever $x=(x_{k})\in c_{0}(\widehat{F},p)$
if and only if $D y\in \ell_{\infty}$ whenever $y=(y_{k})\in
c_{0}(p)$. This means that $a=(a_{n})\in
\{c_{0}(\widehat{F},p)\}^{\gamma}$ whenever $x=(x_{n})\in
c_{0}(\widehat{F},p)$ if and only if $D\in
(c_{0}(p):\ell_{\infty})$. Therefore we obtain Lemma \ref{l2.5} that
$$
\{c_{0}(\widehat{F},p)\}^{\gamma}=\widehat{F}_{3}(p)
$$
and this completes the proof.
\end{pf}
\begin{thm}
Let $K^{*}=\{k\in \mathbb{N}: 0\leq k\leq n\}\cap K$ for $K\in
\mathcal{F}$ and $B\in \mathbb{N}_{2}$. Define the sets
$\widehat{F}_{8}(p),\widehat{F}_{9}(p),\widehat{F}_{10}(p)$ and
$\widehat{F}_{11}(p)$ as follows:
\begin{eqnarray*}
\widehat{F}_{8}(p)&=&\bigcap_{B>1} \bigg\{a=(a_{k})\in \omega:
\sup_{K\in \mathcal{F}} \sum_{n} \bigg|\sum_{k\in
K^{*}}\sum_{j=k}^{n}
\frac{f_{j+1}^{2}}{f_{k}f_{k+1}}a_{j}B^{1/p_{k}}\bigg|<\infty\bigg\}\\
\widehat{F}_{9}(p)&=&\bigcap_{B>1} \bigg\{a=(a_{k})\in
\omega:\sup_{n\in \mathbb{N}} \sum_{k=0}^{n} \bigg|\sum_{j=k}^{n}
\frac{f_{j+1}^{2}}{f_{k}f_{k+1}}a_{j}\bigg|B^{1/p_{k}}<\infty\bigg\}\\
\widehat{F}_{10}(p)&=& \bigcap_{B>1} \bigg\{a=(a_{k})\in \omega:
\exists (\alpha_{k})\subset \mathbb{R}\ni \lim_{n\rightarrow \infty}
\sum_{k=0}^{n} \bigg|\sum_{j=k}^{n}
\frac{f_{j+1}^{2}}{f_{k}f_{k+1}}a_{j}-\alpha_{k}\bigg|B^{1/p_{k}}=0
\bigg\}\\
\widehat{F}_{11}(p)&=&\bigcap_{B>1} \bigg\{a=(a_{k})\in \omega:
\sup_{n\in \mathbb{N}}\sum_{k=0}^{n} \bigg|\sum_{j=k}^{n}
\frac{f_{j+1}^{2}}{f_{k}f_{k+1}}a_{j}\bigg|B^{1/p_{k}}<\infty\bigg\}
\end{eqnarray*}
Then,

(i) $\{\ell_{\infty}(\widehat{F},p)\}^{\alpha}=\widehat{F}_{8}(p)$\\

(ii)
$\{\ell_{\infty}(\widehat{F},p)\}^{\beta}=\widehat{F}_{9}(p)\cap
\widehat{F}_{10}(p)$\\

(iii)
$\{\ell_{\infty}(\widehat{F},p)\}^{\gamma}=\widehat{F}_{11}(p)$.
\end{thm}
\begin{pf}
This may be obtained in the similar way, as mentioned in the proof
of Theorem \ref{t2.9} with Lemmas \ref{l2.6}(i), \ref{l2.7}(i),
\ref{l2.8} instead of Lemmas \ref{l2.3}-\ref{l2.5}. So, we omit the
details.
\end{pf}
\begin{thm}\label{t2.11}
Let $K^{*}=\{k\in \mathbb{N}: 0\leq k\leq n\}\cap K$ for $K\in
\mathcal{F}$ and $B\in \mathbb{N}_{2}$. Define the sets
$\widehat{F}_{12}(p),\widehat{F}_{13}(p),\widehat{F}_{14}(p),\widehat{F}_{15}(p)$
and $\widehat{F}_{16}(p)$ as follows:
\begin{eqnarray*}
\widehat{F}_{12}(p)&=&\bigg\{a=(a_{k})\in \omega: \sup_{K\in
\mathcal{F}} \sup_{k\in \mathbb{N}} \bigg|\sum_{n\in K^{*} }
\sum_{j=k}^{n}
\frac{f_{j+1}^{2}}{f_{k}f_{k+1}}a_{j}\bigg|^{p_{k}}<\infty\bigg\}\\
\widehat{F}_{13}(p)&=&\bigcup_{B>1} \bigg\{a=(a_{k})\in \omega:
\sup_{K\in \mathcal{F}} \sum_{k} \bigg|\sum_{n\in K} \sum_{j=k}^{n}
\frac{f_{j+1}^{2}}{f_{k}f_{k+1}}a_{j}B^{-1}\bigg|^{p_{k}^{'}}<\infty\bigg\}\\
\widehat{F}_{14}(p)&=&\bigcup_{B>1} \bigg\{a=(a_{k})\in \omega:
\sup_{n\in \mathbb{N}} \sum_{k=0}^{n}  \bigg|\sum_{j=k}^{n}
\frac{f_{j+1}^{2}}{f_{k}f_{k+1}}a_{j}B^{-1}\bigg|^{p_{k}^{'}}<\infty\bigg\}\\
\widehat{F}_{15}(p)&=&\bigg\{a=(a_{k})\in \omega: \sup_{n,k\in
\mathbb{N}} \bigg|\sum_{j=k}^{n}
\frac{f_{j+1}^{2}}{f_{k}f_{k+1}}a_{j}\bigg|^{p_{k}}<\infty\bigg\}\\
\widehat{F}_{16}(p)&=&\bigg\{a=(a_{k})\in \omega: \lim_{n\rightarrow
\infty} \sum_{j=k}^{n} \frac{f_{j+1}^{2}}{f_{k}f_{k+1}}a_{j} \quad
\textrm{exists}\bigg\}
\end{eqnarray*}
Then,

(i)
$$
\{\ell (\widehat{F},p)\}^{\alpha}=\left\{\begin{array}{ll}
  \displaystyle \widehat{F}_{12}(p), & 0<p_{k}\leq 1
  \\
    \widehat{F}_{13}(p), & 1<p_{k}\leq H<\infty
\end{array}\right.
$$

(ii)
$$
\{\ell (\widehat{F},p)\}^{\gamma}=\left\{\begin{array}{ll}
  \displaystyle \widehat{F}_{15}(p), & 0<p_{k}\leq 1
  \\
    \widehat{F}_{14}(p), & 1<p_{k}\leq H<\infty.
\end{array}\right.
$$

(iii) Let $0<p_{k}\leq H<\infty$. Then,
$$
\{\ell (\widehat{F},p)\}^{\beta}=\widehat{F}_{14}(p)\cap
\widehat{F}_{15}(p)\cap \widehat{F}_{16}(p).
$$
\end{thm}
\begin{pf}
This may be obtained in the similar way, as mentioned in the proof
of Theorem \ref{t2.9} with Lemmas \ref{l2.6}(ii), \ref{l2.7}(ii),
\ref{l2.8} instead of Lemmas \ref{l2.3}-\ref{l2.5}. So, we omit the
details.
\end{pf}

Now, we may give the sequence of the points of the spaces
$c_{0}(\widehat{F},p),$ $\ell(\widehat{F},p)$ and $c(\widehat{F},p)$
which forms a Schauder basis for those spaces. Because of the
isomorphism $T$, defined in the proof of Theorem \ref{t2.2}, between
the sequence spaces $c_{0}(\widehat{F},p)$ and $c_{0}(p)$,
$\ell(\widehat{F},p)$ and $\ell(p)$, $c(\widehat{F},p)$ and $c(p)$
is onto, the inverse image of the basis of the spaces
$c_{0}(p),\ell(p)$ and $c(p)$ is the basis for our new spaces
$c_{0}(\widehat{F},p),$ $\ell(\widehat{F},p)$ and
$c(\widehat{F},p)$, respectively. Therefore, we have:
\begin{thm}
Let $\mu_{k}=(\widehat{F}x)_k$ for all $k\in \mathbb{N}$. We define
the sequence $b^{(k)}=\{b_{n}^{(k)}\}_{n\in \mathbb{N}}$ for every
fixed $k\in \mathbb{N}$ by
$$
b_{n}^{(k)}=\left\{\begin{array}{ll}
  \displaystyle \frac{f_{n+1}^{2}}{f_{k}f_{k+1}}, & n\geq k,\\
    0, & n<k.
\end{array}\right.
$$
Then,\\
(a) The sequence $\{b^{(k)}\}_{k\in \mathbb{N}}$ is a basis for the
space $c_{0}(\widehat{F},p)$ and any $x\in c_{0}(\widehat{F},p)$ has
a unique representation in the form
$$
x=\sum_{k} \mu_{k} b^{(k)}.
$$
(b) The sequence $\{b^{(k)}\}_{k\in \mathbb{N}}$ is a basis for the
space $\ell(\widehat{F},p)$ and any $x\in \ell(\widehat{F},p)$ has a
unique representation in the form
$$
x=\sum_{k} \mu_{k} b^{(k)}.
$$
(c) The set $\{z,b^{(k)}\}$ is a basis for the space
$c(\widehat{F},p)$ and any $x\in c(\widehat{F},p)$ has a unique
representation in the form
$$
x=lz+\sum_{k} (\mu_{k}-l)b^{(k)}
$$
where $l=\lim_{k\rightarrow \infty} (\widehat{F}x)_k$ and $z=(z_k)$
with
$$
z_{k}=\sum_{j=0}^{k} \frac{f_{k+1}^{2}}{f_{j}f_{j+1}}.
$$
\end{thm}
\section{Some Matrix Mappings on the Sequence Spaces $c_{0}(\widehat{F},p),c(\widehat{F},p)$ $,\ell_{\infty}(\widehat{F},p)$ and $\ell(\widehat{F},p)$ }

In this section, we characterize some matrix mappings on the spaces
$c_{0}(\widehat{F},p),c(\widehat{F},p),\ell_{\infty}(\widehat{F},p)$
and $\ell(\widehat{F},p)$. Firstly, we may give the following
theorem which is useful for deriving the characterization of the
certain matrix classes.
\begin{thm}\cite[Theorem 4.1]{mkfb}
Let $\lambda$ be an FK-space, $U$ be a triangle, $V$ be its inverse
and $\mu$ be arbitrary subset of $\omega$. Then we have $A\in
(\lambda_{U}:\mu)$ if and only if
\begin{equation}
E^{(n)}=(e_{mk}^{(n)})\in (\lambda:c) \quad \textrm{for all} \quad
n\in \mathbb{N}
\end{equation}
and
\begin{equation}
E=(e_{nk})\in (\lambda:\mu)
\end{equation}
where
$$
e_{mk}^{(n)}=\left\{\begin{array}{ll}
  \displaystyle \sum_{j=k}^{m} a_{nj}v_{jk}, & 0\leq k\leq m,\\
    0, & k>m,
\end{array}\right.
$$
and
$$
e_{nk}=\sum_{j=k}^{\infty} a_{nj}v_{jk} \quad \textrm{for all} \quad
k,m,n \in \mathbb{N}.
$$
\end{thm}

Now, we may quote our theorems on the characterization of some
matrix classes concerning with the sequence spaces
$c_{0}(\widehat{F},p),c(\widehat{F},p)$ and
$\ell_{\infty}(\widehat{F},p)$. The necessary and sufficient
conditions characterizing the matrix mappings between the sequence
spaces of Maddox are determined by Grosse-Erdmann \cite{kgge}. Let
$N$ and $K$ denote the finite subset of $\mathbb{N}$, $L$ and $M$
also denote the natural numbers. Prior to giving the theorems, let
us suppose that $(q_{n})$ is a non-decreasing bounded sequence of
positive numbers and consider the following conditions:
\begin{equation}\label{mt23}
\lim_{m\rightarrow \infty} \sum_{j=k}^{m}
\frac{f_{j+1}^{2}}{f_{k}f_{k+1}}a_{nj}=e_{nk},
\end{equation}
\begin{equation}\label{mt24}
\forall L, \quad \sum_{k} |e_{nk}|L^{1/p_{k}}<\infty,
\end{equation}
\begin{equation}\label{mt25}
\exists (\alpha_{k})\subset \mathbb{R}\ni \lim_{m\rightarrow
\infty}\bigg| \sum_{j=k}^{m}
\frac{f_{j+1}^{2}}{f_{k}f_{k+1}}a_{nj}-\alpha_{k}\bigg|=0 \quad
\textrm{for all}\quad  k\in \mathbb{N},
\end{equation}
\begin{equation}\label{mt26}
\exists M, \quad \sup_{m\in \mathbb{N}} \sum_{k=0}^{m}
\bigg|\sum_{j=k}^{m}
\frac{f_{j+1}^{2}}{f_{k}f_{k+1}}a_{nj}\bigg|M^{-1/p_{k}}<\infty,
\end{equation}
\begin{equation}\label{mt27}
\forall L, \exists M, \sup_{m\in \mathbb{N}} \sum_{k=0}^{m}
\bigg|\sum_{j=k}^{m} \frac{f_{j+1}^{2}}{f_{k}f_{k+1}}a_{nj}\bigg|
L^{1/q_{n}} M^{-1/p_{k}}<\infty,
\end{equation}
\begin{equation}\label{mt28}
\lim_{m\rightarrow \infty} \sum_{k} \bigg|\sum_{j=k}^{m}
\frac{f_{j+1}^{2}}{f_{k}f_{k+1}}a_{nj}-\alpha\bigg|=0,
\end{equation}
\begin{equation}\label{mt29}
\forall L, \quad \sup_{n\in \mathbb{N}} \sum_{k}
|e_{nk}|L^{1/p_{k}}<\infty,
\end{equation}
\begin{equation}\label{mt30}
\lim_{n\rightarrow \infty} e_{nk}=\alpha_{k} \qquad \textrm{for all}
\quad k\in \mathbb{N},
\end{equation}
\begin{equation}\label{mt31}
\forall L, \quad \lim_{n\rightarrow \infty} \sum_{k}
|e_{nk}|L^{1/p_{k}}<\infty,
\end{equation}
\begin{equation}\label{mt32}
\forall L, \quad \lim_{n\rightarrow \infty} \sum_{k}
|e_{nk}|L^{1/p_{k}}=0,
\end{equation}
\begin{equation}\label{mt33}
\exists M, \quad \sup_{n\in \mathbb{N}} \bigg(\sum_{k\in K}
|e_{nk}|M^{-1/p_{k}}\bigg)^{q_{n}}<\infty,
\end{equation}
\begin{equation}\label{mt34}
\lim_{n\rightarrow \infty} |e_{nk}|^{q_{n}}=0, \quad \textrm{for
all} \quad k\in \mathbb{N},
\end{equation}
\begin{equation}\label{mt35}
\forall L, \exists M, \quad \sup_{n\in \mathbb{N}} \sum_{k}
|e_{nk}|L^{1/q_{n}}M^{-1/p_{k}}<\infty,
\end{equation}
\begin{equation}\label{mt36}
\lim_{n\rightarrow \infty} |e_{nk}-\alpha_{k}|^{q_{n}}=0, \quad
\textrm{for all} \quad k\in \mathbb{N},
\end{equation}
\begin{equation}\label{mt37}
\exists M, \quad \sup_{n\in \mathbb{N}} \sum_{k}
|e_{nk}|M^{-1/p_{k}}<\infty,
\end{equation}
\begin{equation}\label{mt38}
\forall L, \exists M, \quad \sup_{n\in \mathbb{N}} \sum_{k}
|e_{nk}-\alpha_{k}|L^{1/q_{n}}M^{-1/p_{k}}<\infty,
\end{equation}
\begin{equation}\label{mt39}
\sup_{n\in \mathbb{N}} \bigg|\sum_{k} e_{nk}\bigg|^{q_{n}}<\infty,
\end{equation}
\begin{equation}\label{mt40}
\lim_{n\rightarrow \infty} \bigg|\sum_{k} e_{nk}\bigg|^{q_{n}}=0,
\end{equation}
\begin{equation}\label{mt41}
\lim_{n\rightarrow \infty} \bigg|\sum_{k}
e_{nk}-\alpha\bigg|^{q_{n}}=0,
\end{equation}
\begin{thm}

(i) $A\in (\ell_{\infty}(\widehat{F},p):\ell_{\infty})$ if and only
if (\ref{mt23}),(\ref{mt24}) and (\ref{mt29}) hold.\\

(ii) $A\in (\ell_{\infty}(\widehat{F},p):c)$ if and only if
(\ref{mt23}),(\ref{mt24}), (\ref{mt30}) and (\ref{mt31}) hold.\\

(iii) $A\in (\ell_{\infty}(\widehat{F},p):c_{0})$ if and only if
(\ref{mt23}),(\ref{mt24}) and (\ref{mt32}) hold.
\end{thm}
\begin{thm}

(i) $A\in (c_{0}(\widehat{F},p):\ell_{\infty}(q))$ if and only if
(\ref{mt25}), (\ref{mt26}), (\ref{mt27}) and (\ref{mt33}) hold.\\

(ii)  $A\in (c_{0}(\widehat{F},p):c_{0}(q))$ if and only if
(\ref{mt25}), (\ref{mt26}), (\ref{mt27}), (\ref{mt34}) and
(\ref{mt35}) hold.\\

(iii) $A\in (c_{0}(\widehat{F},p):c(q))$ if and only if
(\ref{mt25}), (\ref{mt26}), (\ref{mt27}), (\ref{mt36}), (\ref{mt37})
and (\ref{mt38}) hold.

\end{thm}
\begin{thm}
(i)  $A\in (c(\widehat{F},p):\ell_{\infty}(q))$ if and only if
(\ref{mt25}), (\ref{mt26}), (\ref{mt27}), (\ref{mt28}), (\ref{mt33})
and (\ref{mt39}) hold.\\

(ii) $A\in (c(\widehat{F},p):c_{0}(q))$ if and only if (\ref{mt25}),
(\ref{mt26}), (\ref{mt27}), (\ref{mt28}), (\ref{mt34}), (\ref{mt35})
and (\ref{mt40}) hold.\\

(iii) $A\in (c(\widehat{F},p):c(q))$ if and only if (\ref{mt25}),
(\ref{mt26}), (\ref{mt27}), (\ref{mt28}), (\ref{mt36}),
(\ref{mt37}), (\ref{mt38}) and (\ref{mt41}) hold.
\end{thm}




\end{document}